\newcommand{\al}{\alpha}               
\newcommand{\be}{\beta}
\newcommand{\ga}{\gamma}               

\newcommand{\de}{\delta}               

\newcommand{\lb}{\lambda}              

\newcommand{\sig}{\sigma}              

\newcommand{\om}{\omega}               

\newcommand{\veps}{\varepsilon}        
\newcommand{\vphi}{\varphi}

\newcommand{\cal}{\mathcal}

\newcommand{\calf}{{\cal F}}

\newcommand{\Dom}{{\rm Dom}}
\newcommand{\Fix}{{\rm Fix}}

\newcommand{\es}{\emptyset}

\newcommand{\limpl}{\Longrightarrow}
    
\newcommand{\lequi}{\Longleftrightarrow}
\newcommand{\oo}{\infty}

\def\dtends   {\stackrel {\it d}{\longrightarrow}}

\newcommand{\barr}{\begin{array}}        
\newcommand{\earr}{\end{array}}
\newcommand{\bcor}{\begin{corollary}}    
\newcommand{\ecor}{\end{corollary}}
\newcommand{\ben}{\begin{enumerate}}     
\newcommand{\een}{\end{enumerate}}
\newcommand{\beq}{\begin{equation}}      
\newcommand{\eeq}{\end{equation}}
\newcommand{\bex}{\begin{example}}     
\newcommand{\eex}{\end{example}}
\newcommand{\bit}{\begin{itemize}}       
\newcommand{\eit}{\end{itemize}}
\newcommand{\blemma}{\begin{lemma}}      
\newcommand{\elemma}{\end{lemma}}
\newcommand{\bproof}{\begin{proof}}      
\newcommand{\eproof}{\end{proof}}
\newcommand{\bprop}{\begin{proposition}} 
\newcommand{\eprop}{\end{proposition}}
\newcommand{\brem}{\begin{remark}} 
\newcommand{\erem}{\end{remark}}
\newcommand{\btab}{\begin{tabular}}      
\newcommand{\etab}{\end{tabular}}
\newcommand{\btheorem}{\begin{theorem}}  
\newcommand{\etheorem}{\end{theorem}}

\documentclass[reqno]{amsart}

\newtheorem{theorem}{\bf Theorem}
\newtheorem{corollary}{\bf Corollary}
\newtheorem{example}{\bf Example}
\newtheorem{lemma}{\bf Lemma}
\newtheorem{proposition}{\bf Proposition}
\newtheorem{remark}{\bf Remark}

\begin{document}

\title
[Weakly Contractive Maps in Altering Metric Spaces]
{WEAKLY CONTRACTIVE MAPS IN \\
ALTERING METRIC SPACES}

\author{Mihai Turinici}
\address{
"A. Myller" Mathematical Seminar;
"A. I. Cuza" University;
700506 Ia\c{s}i, Romania
}
\email{mturi@uaic.ro}

% \date{7 October, 2013}



\begin{abstract}
The weakly contractive metric type 
fixed point result in Berinde
[Nonlinear Anal. Forum, 9 (2004), 45-53]
is "almost" covered by
the  related altering metric one due to 
Khan et al
[Bull. Austral. Math. Soc., 30 (1984), 1-9].
Further extensions of these statements 
are then provided.
\end{abstract}

\maketitle

%%%%%%%%%%%%%%%%%%%%%%%%%%%%%%%%%%%%%%%

% Section 1
\section{Introduction}
\setcounter{equation}{0}

Let $X$ be a nonempty set.
Call the subset $Y$ of $X$, 
{\it almost singleton} (in short: {\it asingleton})
provided [$y_1,y_2\in Y$ implies $y_1=y_2$];
and {\it singleton},
if, in addition, $Y$ is nonempty;
note that, in this case,
$Y=\{y\}$, for some $y\in X$. 
Further, let $d:X\times X\to R_+:=[0,\oo[$ 
be a {\it metric} over it;
the couple $(X,d)$ 
will be termed a {\it metric space}.
Finally, let $T\in \calf(X)$ be a selfmap of $X$.
[Here, for each couple $A, B$ of nonempty sets,
$\calf(A,B)$ stands for the class of all functions 
from $A$ to $B$; 
when $A=B$, we write $\calf(A)$ in place of $\calf(A,A)$].
Denote $\Fix(T)=\{x\in X; x=Tx\}$;
each point of this set is referred to as 
{\it fixed} under $T$.
In the metrical fixed point theory, 
such points are to be determined 
according to the context below,
comparable with the one in 
Rus \cite[Ch 2, Sect 2.2]{rus-2001}:

{\bf 1a)} 
We say that $T$ is a {\it Picard operator} (modulo $d$) if,
for each $x\in X$, 
the iterative sequence
$(T^nx; n\ge 0)$ is $d$-convergent

{\bf 1b)} 
We say that $T$ is a {\it strong Picard operator} (modulo $d$) if,
for each $x\in X$, $(T^nx; n\ge 0)$ is $d$-convergent 
and $\lim_n(T^nx)$ belongs to $\Fix(T)$

{\bf 1c)}
We say that $T$ is a {\it globally strong Picard operator} (modulo $d$) if
it is a strong Picard operator (modulo $d$),
and
$\Fix(T)$ is an asingleton (hence, a singleton).

In this perspective, 
a basic answer to the posed question is the 1922 one due to 
Banach \cite{banach-1922}:
it states that, 
if  $T$ is {\it $(d;\al)$-contractive}, i.e.,
\ben
\item[] (a01)\ \ 
$d(Tx,Ty)\le \al d(x,y)$,\ \ for all  $x,y\in X$,
\een
for some $\al\in [0,1[$, 
and $(X,d)$ is {\it complete},
then $T$ is a globally strong Picard operator (modulo $d$).

This result found a multitude of applications in 
operator equations theory;
so, it was the subject of many extensions.
For example, a natural way of doing this  
is by considering (implicit) "functional" contractive conditions
of the form
\ben
\item[] (a02)\ \ 
$F(d(Tx,Ty),d(x,y),d(x,Tx),d(y,Ty),d(x,Ty),d(y,Tx))\le 0$,\\
for all  $x,y\in X$;
\een
where $F:R_+^6\to R$ is an appropriate function.
For more details about the possible choices of $F$
we refer to the 1977 paper by
Rhoades \cite{rhoades-1977};
see also 
Turinici \cite{turinici-1976}.
Here, we shall be concerned with a 2004 contribution 
in the area due to Berinde \cite{berinde-2004}.
Given $\al,\lb\ge 0$, 
let us say that  $T$ is a {\it weak $(\al,\lb)$-contraction} (modulo $d$)
provided
\ben
\item[] (a03)\ \ 
$d(Tx,Ty)\le \al d(x,y)+\lb d(Tx,y)$,\ \ for all $x,y\in X$.
\een

% Theorem 1
\btheorem \label{t1}
Suppose that $T$ is a weak $(\al,\lb)$-contraction (modulo $d$),
where $\al\in [0,1[$. 
In addition, let $(X,d)$ be complete.
Then, $T$ is a strong Picard operator (modulo $d$).
\etheorem

In a subsequent paper devoted to the same question,
Berinde \cite{berinde-2003} claims that this class of contractions 
introduced by him is for the first time
considered in the literature.
Unfortunately, his assertion is not true:
conclusions of Theorem \ref{t1} 
are "almost" covered by 
a related 1984 statement due to  
Khan et al \cite{khan-swaleh-sessa-1984},
in the context of altering distances.
This, among others, motivated us to propose
an appropriate extension of the quoted 
statement; details are given in Section 3.
The preliminary material for our device is listed in 
Section 2.
Finally, in Section 4, a "functional" extension of
Berinde's result is established.
Further aspects will be delineated 
in a separate paper.

% Section 2
\section{Preliminaries}
\setcounter{equation}{0}

Let $(X,d)$ be a metric space.
Let us say that the sequence $(x_n)$ in $X$, {\it $d$-converges}
to $x\in X$ (and write: $x_n\dtends x$) iff $d(x_n,x)\to 0$; that is
\ben
\item[] (b01)\ \ 
$\forall \veps> 0$, $\exists p=p(\veps)$, $\forall n$:\ \ 
($p\le n$ $\limpl$ $d(x_n,x)\le \veps$).
\een
The subset $\lim_n(x_n)$  of all such $x$ is an asingleton, 
because $d$ is sufficient;
when it is nonempty, 
$(x_n)$ is called {\it $d$-convergent}.
Note that, in this case,  $\lim_n(x_n)$ is a singleton,
$\{z\}$; as usually, we write $\lim_n(x_n)=z$.
Further, let us say that $(x_n)$ is 
{\it $d$-Cauchy}, provided 
$d(x_m,x_n)\to 0$ as $m,n\to \oo$, $m< n$; that is
\ben
\item[] (b02)\ \ 
$\forall \veps> 0$, $\exists q=q(\veps)$, $\forall (m,n)$:\ \ 
($q\le m< n$ $\limpl$ $d(x_m,x_n)\le \veps$).
\een
Clearly, any $d$-convergent sequence is $d$-Cauchy too;
when the reciprocal holds as well, 
$(X,d)$ is called {\it complete}.
Concerning this aspect, 
note that any $d$-Cauchy sequence $(x_n; n\ge 0)$ 
is {\it $d$-semi-Cauchy}; i.e.,
\ben
\item[] (b03)\ \ 
$\rho_n:=d(x_n,x_{n+1})\to 0$,\  as $n\to \oo$.
\een
The following result involving this property is useful in the sequel.
For each sequence $(r_n; n\ge 0)$ in $R$ and each $r\in R$, put
\ben
\item[]
$r_n \to r+$ iff [$r_n> r$, $\forall n$] and $r_n\to r$.
\een

% Proposition 1
\bprop \label{p1}
Suppose that $(x_n; n\ge 0)$ is $d$-semi-Cauchy, 
but not $d$-Cauchy.
There exists then $\eta> 0$,  $j(\eta)\ge 0$,
and a couple of rank sequences
$(m(j); j\ge 0)$, $(n(j); j\ge 0)$,
in such a way that
\beq \label{201}
j\le m(j)< n(j),\ \  \al(j):=d(x_{m(j)},x_{n(j)})> \eta,\ \forall j\ge 0
\eeq
\beq \label{202}
n(j)-m(j)\ge 2,\ \be(j):=d(x_{m(j)},x_{n(j)-1})\le \eta,\
\forall j\ge j(\eta)
\eeq
\beq \label{203}
\mbox{
$\al(j)\to \eta+$\ 
(hence, $\al(j)\to \eta$),\ as $j\to \oo$
}
\eeq
\beq \label{204}
\al_{p,q}(j):=d(x_{m(j)+p},x_{n(j)+q})\to \eta,\ \mbox{as}\ 
j\to \oo,\ \forall p,q\in \{0,1\}.
\eeq
\eprop

A  proof of this may be found in  
Khan et al \cite{khan-swaleh-sessa-1984}.
For completeness reasons, we supply 
an argument which differs, in part, from the original one.

\bproof {\bf (Proposition \ref{p1})}
As (b02) does not hold, there exists $\eta> 0$ with
$$
A(j):=\{(m,n)\in N\times N; j\le m< n, d(x_m,x_n)> \eta\}
\ne \es,\ \forall j\ge 0.
$$
Having this precise, denote, for each $j\ge 0$,
\ben
\item[]
$m(j)=\min \Dom(A(j))$,\ $n(j)=\min A(m(j))$.
\een
As a consequence, the couple of rank-sequences 
$(m(j); j\ge 0)$, $(n(j); j\ge 0)$  fulfills (\ref{201}).
On the other hand, letting 
the index $j(\eta)\ge 0$ be such that 
\beq \label{205}
(\rho_k=)\ d(x_k,x_{k+1})< \eta,\ \  \forall k\ge j(\eta),
\eeq
it is clear that (\ref{202}) holds too.
Finally, by the triangular property,
$$ 
\eta < \al(j)\le
\be(j)+\rho_{n(j)-1}\le \eta+ \rho_{n(j)-1},\ \ \forall j\ge j(\eta);
$$
and this yields (\ref{203}); hence, the case 
$(p=0,q=0)$ of (\ref{204}).
Combining with 
$$ 
\al(j)-\rho_{n(j)}\le d(x_{m(j)},x_{n(j)+1}) 
\le \al(j)+\rho_{n(j)},\ \forall j\ge j(\eta)
$$
establishes the case $(p=0,q=1)$ of the same. 
The remaining situations are deductible in a similar way.
\eproof

% Section 3
\section{Main result}
\setcounter{equation}{0}

Let $X$ be a nonempty set; 
and $d(.,.)$ be a metric over it [in the usual sense].
Further, let $\vphi\in \calf(R_+)$ be an 
{\it altering function}; i.e.
\ben
\item[] (c01)\ \ 
$\vphi$ is continuous, increasing, 
and reflexive sufficient [$\vphi(t)=0$ iff $t=0$].
\een
The associated map (from $X\times X$ to $R_+$)
\ben
\item[] (c02)\ \ 
$e(x,y)=\vphi(d(x,y))$,\ \ $x,y\in X$
\een
has the  immediate properties
\beq \label{301}
\mbox{
$e(x,y)=e(y,x),\ \forall x,y\in X$  \quad ($e$ is symmetric)
}
\eeq
\beq \label{302}
\mbox{
$e(x,y)=0$ $\lequi$ $x=y$  \quad ($e$ is reflexive sufficient). 
}
\eeq
So, it is a (reflexive sufficient) {\it symmetric}, under 
the Hicks-Rhoades terminology \cite{hicks-rhoades-1999}.
In general, $e(.,.)$ is not endowed with the triangular property; 
but, in compensation to this, one has 
(as $\vphi$ is increasing and continuous)
\beq \label{303}
e(x,y)> e(u,v) \limpl d(x,y)> d(u,v)
\eeq
\beq \label{304}
\mbox{
$x_n\dtends x$, $y_n\dtends y$\ \ \mbox{implies}\ \ 
$e(x_n,y_n)\to e(x,y)$.
}
\eeq

Let $T\in \calf(X)$ be a selfmap of $X$.
The problem
involving its fixed points is the already stated one.
In the following, we are trying to solve it, under the
precise "altering" context. Denote, for $x,y\in X$, 
\ben
\item[] (c03)\ \  
$M_1(x,y)=e(x,y)$,\
$M_2(x,y)=(1/2)[e(x,Tx)+e(y,Ty)]$,\\
$M_3(x,y)=\min\{e(x,Ty),e(Tx,y)\}$,\\
$M(x,y)=\max\{M_1(x,y),M_2(x,y),M_3(x,y)\}$.
\een
Further, given $\psi\in \calf(R_+)$,  
we say that $T$ is {\it $(d,e;M;\psi)$-contractive},
provided
\ben
\item[] (c04)\ \  
$e(Tx,Ty)\le \psi(d(x,y))M(x,y)$,\ 
$\forall x,y \in X$, $x\ne y$.
\een
The properties of $\psi$ to be used here are the following
\ben
\item[] (c05)\ \ 
$\psi$ is strictly subunitary on $R_+^0:=]0,\oo[$:\
$\psi(s)< 1$, $\forall s\in R_+^0$
\item[] (c06)\ \ 
$\psi$ is right Boyd-Wong on $R_+^0$:\
$\limsup_{t\to s+} \psi(t)< 1$,\ 
$\forall s\in R_+^0$.
\een
These are related to the developments of 
Boyd and Wong \cite{boyd-wong-1969};
we do not give details.

Our main result in this exposition is

% Theorem 2
\btheorem \label{t2}
Suppose that $T$ is $(d,e;M;\psi)$-contractive,
where $\psi\in \calf(R_+)$ is 
strictly subunitary and right Boyd-Wong on $R_+^0$.
In addition, let $(X,d)$ be complete.
Then, $T$ is a globally strong Picard operator (modulo $d$).
\etheorem

\bproof
First, we check 
the asingleton property for $\Fix(T)$.
Let $z_1,z_2\in \Fix(T)$ be such that $z_1\ne z_2$;
hence $\de:=d(z_1,z_2)> 0$, $\veps:=e(z_1,z_2)> 0$.
By definition, 
$$
M_1(z_1,z_2)=\veps,\
M_2(z_2,z_2)=0,\ M_3(x,y)=\veps;\  
\mbox{hence}\
M(x,y)=\veps.
$$
By the contractive condition (written at $(z_1,z_2)$)
$$
\veps=e(z_1,z_2)=e(Tz_1,Tz_2)\le 
\psi(\de)M(z_1,z_2)=\psi(\de)\veps;
$$
hence, $1\le \psi(\de)< 1$; contradiction;
and the asingleton property follows.
It remains now to verify the strong Picard property
(modulo $d$).
Fix a certain $x_0\in X$; and put $(x_n=T^nx_0; n\ge 0)$.
If $x_n=x_{n+1}$ for some $n\ge 0$, we are done;
so, without loss, one may assume that
\ben
\item[] (c07)\ \ 
$\rho_n:=d(x_n,x_{n+1})> 0$\ 
(hence, $\sig_n:=e(x_n,x_{n+1})> 0$),\ for all $n$.
\een
There are several steps to be passed.

{\bf I)}
For the arbitrary fixed $n\ge 0$, we have
$$  \barr{l}
M_1(x_n,x_{n+1})=\sig_n, \\
M_2(x_n,x_{n+1})=(1/2)[\sig_n+\sig_{n+1}]\le 
\max\{\sig_n,\sig_{n+1}\}, \\
M_3(x_n,x_{n+1})=0;\ 
\mbox{hence}\ 
M(x_n,x_{n+1})\le \max\{\sig_n,\sig_{n+1}\}.
\earr
$$
By the contractive condition (written at $(x_n,x_{n+1})$), 
$$
\sig_{n+1}\le \psi(\rho_n)\max\{\sig_n,\sig_{n+1}\},\
\forall n.
$$
This, along with (c07),  yields 
(as $\psi$ is strictly subunitary on $R_+^0$) 
\beq \label{305}
\sig_{n+1}/\sig_n \le \psi(\rho_n)< 1,\ \forall n.
\eeq
As a direct consequence,
$$
\mbox{
$\sig_n> \sig_{n+1}$\ 
(hence, $\rho_n> \rho_{n+1}$),\ for all $n$.
}
$$
The sequence
$(\rho_n; n\ge 0)$ is therefore 
strictly descending in $R_+$; hence, 
$\rho:=\lim_n(\rho_n)$ exist in $R_+$ and
$\rho_n> \rho$, $\forall n$.
Likewise, the sequence 
$(\sig_n=\vphi(\rho_n); n\ge 0)$ is
strictly descending in $R_+$; hence, 
$\sig:=\lim_n(\sig_n)$ exists;
with, in addition,  $\sig=\vphi(\rho)$. 
We claim that $\rho=0$. 
Assume by contradiction that $\rho> 0$; hence $\sig> 0$.
Passing to $\limsup$ as $n\to \oo$ in (\ref{305}), yields
$$
1\le \limsup_n \psi(\rho_n)\le \limsup_{t\to \rho+}\psi(t)< 1;
$$
contradiction. Hence, $\rho=0$; i.e.,
\beq \label{306}
\rho_n:=d(x_n,x_{n+1})\to 0,\ \mbox{as}\ n\to \oo.
\eeq

{\bf II)}
We now show that $(x_n; n\ge 0)$ is $d$-Cauchy.
Suppose that this is not true.
By Proposition \ref{p1}, 
there exist $\eta> 0$,  $j(\eta)\ge 0$,
and a couple of rank sequences
$(m(j); j\ge 0)$, $(n(j); j\ge 0)$,
in such a way that (\ref{201})-(\ref{204}) hold.
Denote for simplicity $\zeta=\vphi(\eta)$; 
hence, $\zeta> 0$.
By the notations used there,
we may write as $j\to \oo$
$$
\lb_j:=e(x_{m(j)+1},x_{n(j)+1})=\vphi(\al_{1,1}(j))\to \zeta. 
$$
In addition, we have (again under $j\to \oo$)
$$
\barr{l}
M_1(x_{m(j)},x_{n(j)})=\vphi(\al(j))\to \zeta \\
M_2(x_{m(j)},x_{n(j)})=
(1/2)[\vphi(\rho_{m(j)})+\vphi(\rho_{n(j)})]\to 0 \\
M_3(x_{m(j)},x_{n(j)})=
\min\{\vphi(\al_{0,1}(j)),\vphi(\al_{1,0}(j))\}\to \zeta; 
\earr
$$
and this, by definition, yields
$$
(0<)\ \mu_j:=M(x_{m(j)},x_{n(j)})\to \zeta\ 
\mbox{as}\ j\to \oo.
$$
From the contractive condition
(written at $(x_{m(j)},x_{n(j)})$)
$$
\lb_j/\mu_j\le \psi(\al(j)),\ \forall j\ge j(\eta);
$$
so that, passing to $\limsup$ as $j\to \oo$
$$
1\le \limsup_j\psi(\al(j))\le \limsup_{t\to \eta+} \psi(t)< 1;
$$
contradiction. Hence, $(x_n; n\ge 0)$ is $d$-Cauchy, as claimed.

{\bf III)}
As $(X,d)$ is complete, there exists a (uniquely determined)
$z\in X$ with $x_n\dtends z$; 
hence $\ga_n:=d(x_n,z)\to 0$ as $n\to \oo$.

Two assumptions are open before us:

{\bf i)} For each $h\ge 0$, there exists $k> h$ with $x_k=z$.
In this case, there exists a sequence of ranks $(m(i); i\ge 0)$
with $m(i)\to \oo$  as $i\to \oo$  such  that $x_{m(i)}=z$ 
(hence,  $x_{m(i)+1}=Tz$), for all $i\ge 0$. 
Letting $i$ tends to infinity and using the fact that 
$(y_i:=x_{m(i)+1}; i\ge 0)$  is a subsequence of $(x_i; i\ge 0)$, 
we get $z=Tz$.

{\bf ii)} There exists $h\ge 0$ such that $n\ge h$ $\limpl$ $x_n\ne z$;
hence, $\ga_n> 0$. 
Suppose that $z\ne Tz$; i.e., $\theta:=d(z,Tz)> 0$;
hence, $\om:=e(z,Tz)> 0$.
Note that, in such a case, $\de_n:=d(x_n,Tz)\to \theta$.
From our previous notations, we have (as $n\to \oo$)
$$
\lb_n:=e(x_{n+1},Tz)=\vphi(\de_{n+1})\to \vphi(\theta)=\om.
$$
In addition (again under $n\to \oo$), 
$$
\barr{l}
M_1(x_n,z)=\vphi(\ga_n)\to 0, \
M_2(x_n,z)=(1/2)[\sig_n+\om]\to \om/2 \\
M_3(x_n,z)=\min\{\vphi(\de_{n}),\vphi(\ga_{n+1})\}\to 0;
\earr
$$
wherefrom (as $M_1(x_n,z), M_2(x_n,z)> 0$, $\forall n\ge h$)
$$
\mbox{
$(0<) \mu_n:=M(x_n,z)\to \om/2$,\ as $n\to \oo$.
}
$$
By the contractive condition (written at $(x_n,z)$)
$$
\lb_n\le \psi(\ga_n) \mu_n< \mu_n,\ \forall n\ge h
$$
we then have (passing to limit as $n\to \oo$),
$\om\le \om/2$; hence $\om=0$.
This yields $\theta=0$; contradiction.
Hence, $z$ is fixed under $T$ and 
the proof is complete.
\eproof

In particular, the right Boyd-Wong on $R_+^0$ property
of $\psi$ is assured when this function 
fulfills (c05) and is decreasing on $R_+^0$.
As a consequence, the following particular
version of our main result may be stated.

% Theorem 3
\btheorem \label{t3}
Suppose that $T$ is $(d,e;M;\psi)$-contractive,
where $\psi\in \calf(R_+)$ is 
strictly subunitary and decreasing on $R_+^0$.
In addition, let $(X,d)$ be complete.
Then, $T$ is a globally strong Picard operator (modulo $d$).
\etheorem

Let $a,b,c\in \calf(R_+)$ be a triple of functions.
We say that the selfmap $T$  of $X$ is 
{\it $(d,e;a,b,c)$-contractive} if
\ben
\item[] (c08)\ \ 
$e(Tx,Ty)\le  
a(d(x,y))e(x,y)+ 
b(d(x,y))[e(x,Tx)+e(y,Ty)]+$\\
$c(d(x,y))\min\{e(x,Ty),e(Tx,y)\}$,\ \ 
$\forall x,y\in X$, $x\ne y$.
\een
Denote for simplicity $\psi=a+2b+c$;
it is clear that, under such a condition,
$T$ is $(d,e;M;\psi)$-contractive.
Consequently, the following statement is
a particular case of Theorem \ref{t2} above:

% Theorem 4
\btheorem \label{t4}
Suppose that $T$ is $(d,e;a,b,c)$-contractive,
where the triple of functions $a,b,c\in \calf(R_+)$ 
is such that their associated function $\psi=a+2b+c$
is strictly subunitary and right Boyd-Wong on $R_+^0$.
In addition, let $(X,d)$ be complete.
Then, conclusions of Theorem \ref{t2} hold.
\etheorem

In particular, when $a,b,c$ are all 
decreasing on $R_+^0$, the right Boyd-Wong property 
on $R_+^0$ of $\psi$  holds;
note that, in this case, Theorem \ref{t4}
is also reducible to Theorem \ref{t3}.
This is just the 1984 fixed point result in 
Khan et al \cite{khan-swaleh-sessa-1984}.

Finally, it is worth mentioning 
that the nice contributions of these authors
were the starting point for a series of results involving 
altering contractions, like the ones in
Bhaumik et al \cite{bhaumik-das-metiya-choudhury-2012},
Nashine and Samet \cite{nashine-samet-2011},
or
Sastry and Babu \cite{sastry-babu-1999};
see also
Pathak and Shahzad \cite{pathak-shahzad-2009}.
However, according to the developments in 
Jachymski \cite{jachymski-2011},
most of these
(including the 
Dutta-Choudhury's contribution \cite{dutta-choudhury-2008})
are in fact reducible to standard techniques;
we do not give details.

% Section 4
\section{Further aspects}
\setcounter{equation}{0}

Let again $(X,d)$ be a metric space 
and $T\in \calf(X)$ be a selfmap of $X$.
A basic particular case of Theorem \ref{t4} 
corresponds to the choices 
$\vphi$=identity and [$a,b,c$=constants].
The corresponding form of Theorem \ref{t4} 
is comparable with Theorem \ref{t1}.
However, the inclusion between these is 
not complete.
This raises the question of determining
proper extensions of Theorem \ref{t1},
close enough to Theorem \ref{t4}.
A direct answer to this is provided by

% Theorem 5
\btheorem \label{t5}
Let the numbers $a,b\in R_+$ 
and the function $K\in \calf(R_+)$ 
be such that 
\ben
\item[] (d01)\ \ 
$d(Tx,Ty)\le ad(x,y)+b[d(x,Tx)+d(y,Ty)]+ 
K(d(Tx,y))$,\ $\forall x,y\in X$
\item[] (d02)\ \ 
$a+2b< 1$ and $K(t)\to 0=K(0)$, as $t\to 0$.
\een
In addition, let $(X,d)$ be complete.
Then, $T$ is strong Picard (modulo $d$).
\etheorem

\bproof
Take an arbitrary fixed $u\in X$.
By the very contractive condition (written at $(T^nu,T^{n+1}u)$),
we have the evaluation
\beq \label{401}
d(T^{n+1}u,T^{n+2}u)\le \lb d(T^{n}u,T^{n+1}u),\ \ \forall n\ge 0;
\eeq
where $\lb:=(a+b)/(1-b)< 1$.
This yields 
\beq \label{402}
d(T^nu,T^{n+1}u)\le \lb^n d(u,Tu),\ \ \forall n\ge 0.
\eeq
Consequently, $(T^nu; n\ge 0)$ is $d$-Cauchy; 
whence (by completeness) 
$$
\mbox{
$T^nu\dtends z:=T^\oo u$, for  some $z\in X$.
}
$$
From the contractive condition (written at $(T^nu,z)$),
$$
d(T^{n+1}u,Tz)\le 
ad(T^nu,z)+b[d(T^nu,T^{n+1}u)+d(z,Tz)]+K(d(T^{n+1}u,z)),\ 
\forall n.
$$
Passing to limit as $n\to \oo$ gives (via (d02))
$d(z,Tz)\le bd(z,Tz)$; so that, if $z\ne Tz$,
one gets $1\le b< 1/2$, contradiction.
Hence $z=Tz$; and the proof is complete.
\eproof

In particular, when $b=0$ and $K(.)$ is linear 
($K(t)=\lb t$, $t\in R_+$, for some $\lb\ge 0$),
this result is just Theorem \ref{t1}.
Note that, from (\ref{402}), one has 
for these "limit" fixed points, 
the error approximation formula
\beq \label{403}
d(T^nu,T^\oo u)\le [\lb^n/(1-\lb)]d(u,Tu),\ \ 
\forall u\in X,\ \forall n.
\eeq
However, when $\Fix(T)$ is non-singleton, 
this "local" evaluation is 
without practical effect,
by the highly unstable character of the map 
$u\mapsto T^\oo u$.
In fact, assume for simplicity that 
$T$ is continuous; and fix in the following 
$u_0\in X$. 
Given $\veps> 0$, there exists $\de> 0$, such that 
$x\in X(u_0,\de)$ $\limpl$ $Tx\in X(Tu_0,\veps)$;
here, for each $x\in X$, $\rho> 0$, we denoted 
$X(x,\rho)=\{y\in X; d(x,y)< \rho\}$
(the {\it open sphere} with center $x$ and radius $\rho$).
By (\ref{403}), one gets
a "local-global" relation like
\beq \label{404}
d(T^nu,T^\oo u)\le [\lb^n/(1-\lb)]\mu(u_0),\ 
\forall u\in X(u_0,\de),\ \forall n;
\eeq
where, by definition, 
$\mu(u_0)=\sup\{d(x,Tx); x\in X(u_0,\de)\}$.
Now, in practice, the starting point $u_0$ is 
approximated by a certain $v_0\in X(u_0,\de)$;
with, in general, $v_0\ne u_0$.
Suppose that the iterates 
$(T^n v_0; n\ge 0)$ are 
calculated  in a complete (and exact) way.
The approximation formula (\ref{404}) 
gives, for the point in question,
\beq \label{405}
d(T^nv_0,T^\oo v_0)\le [\lb^n/(1-\lb)]\mu(u_0),\ 
\ \mbox{for all}\ n.
\eeq 
This yields a good evaluation 
for the fixed point $T^\oo v_0$;
but, it may 
have no impact upon 
the fixed point $T^\oo u_0$ 
(that we want to approximate),
as long as it is distinct from 
the preceding fixed point.
On the contrary, when 
$\Fix(T)$ is a singleton, $\{z\}$,
the local-global relation (\ref{404}) becomes
\beq \label{406}
d(T^nu,z)\le [\lb^n/(1-\lb)]\mu(u_0),\ 
\forall u\in X(u_0,\de),\ \forall n.
\eeq
In this case, 
for the approximation $v_0\in X(u_0,\de)$
of the starting point $u_0$, we have,
by (\ref{405}) above,
\beq \label{407}
d(T^nv_0,z)\le [\lb^n/(1-\lb)]\mu(u_0),\ 
\ \mbox{for all}\ n;
\eeq 
hence, the unique fixed point $z$ 
in $\Fix(T)$ is very well evaluated 
by the iterates $(T^n v_0; n\ge 0)$, 
with any degree of accuracy.
Summing up, any such contraction
$T$ is
Hyers-Ulam unstable, whenever
$\Fix(T)$ is not a singleton,
and Hyers-Ulam stable, 
provided $\Fix(T)$ is a singleton.
Some related facts 
may be found in the 1998 monograph by
Hyers, Isac and Rassias \cite{hyers-isac-rassias-1998}.

Note finally that, further enlargements 
of this result are possible,
in the realm of
partial metric spaces,
introduced under the lines in 
Matthews \cite{matthews-1994}. 
Likewise, an extension of these facts
is possible 
to the framework of quasi-ordered metric spaces
under the lines in 
Agarwal et al \cite{agarwal-el-gebeily-o-regan-2008};
see also
Turinici \cite{turinici-1986}.
A development of both these directions 
will be given in a separate paper.

%  References

\end{document}